\documentclass[11pt,a4paper]{article}

\usepackage{amssymb} 
\usepackage{amsfonts}
\usepackage{epsfig}

\def\beq{\begin{equation}}
\def\eeq{\end{equation}}
\def\bea{\begin{eqnarray}}
\def\eea{\end{eqnarray}}
\def\beas{\begin{eqnarray*}}
\def\eeas{\end{eqnarray*}}

\makeindex             

\usepackage{graphics,epsfig,theorem,latexsym,amssymb,amsmath}
\usepackage{times,epic,eepic}

\newenvironment{keywords}{ \noindent {\small\bf Key Words}:}{ }
\newtheorem{theorem}{Theorem}[section]

\begin{document}

\title{A One-Dimensional Local Tuning Algorithm for Solving GO Problems
with Partially Defined Constraints\thanks{This research was
supported by the following grants: FIRB RBNE01WBBB, FIRB
RBAU01JYPN, PRIN 2005017083-002, and RFBR 04-01-00455-a.   The
authors would like to thank anonymous referees for their subtle
suggestions.}}


\newcommand{\nms}{\normalsize}

\author{Yaroslav D.
Sergeyev $^{1,2}$\footnote{Corresponding author} \and Dmitri E.
Kvasov$^{1,2}$ \and Falah M.H.~Khalaf$^{3}$ \\
 $^{1}$DEIS, University of Calabria, Via P. Bucci, 42C\\ 87036
-- Rende (CS), Italy, \\
$^{2}$ Software Department, N.I.~Lobatchevsky State University\\
Nizhni Novgorod, Russia, \\
$^{3}$ Department of Mathematics, University of Calabria, Italy,
\\ \texttt{yaro@si.deis.unical.it}
\\ \texttt{kvadim@si.deis.unical.it}
\\ \texttt{falah@mat.unical.it} }

%
%
\date{}
\maketitle

\begin{abstract}
Lipschitz one-dimensional constrained global optimization (GO)
problems where both the objective function and constraints can be
multiextremal and non-differentiable are considered in this paper.
Problems, where the constraints are verified in an a priori given
order fixed by the nature of the problem are studied. Moreover, if
a constraint is not satisfied at a point, then the remaining
constraints and the objective function can be undefined at this
point. The constrained problem is reduced to a discontinuous
unconstrained problem by the index scheme without introducing
additional parameters or variables. A new geometric method using
adaptive estimates of local Lipschitz constants is introduced. The
estimates are calculated by using the local tuning technique
proposed recently. Numerical experiments show quite a satisfactory
performance of the new method in comparison with the penalty
approach and a method using a priori given Lipschitz constants.
 \end{abstract}

\section{Introduction}

It happens often in engineering optimization problems (see
\cite{Pinter (1996),Strongin (1978),Strongin and Sergeyev (2000)})
that the objective function and constraints can be multiextremal,
non-differentiable, and partially defined. The latter means that
the constraints are verified in a priori given order fixed by the
nature of the problem and if a constraint is not satisfied at a
point, then the remaining constraints and the objective function
can be undefined at this point. This kind of problems is difficult
to solve even in the one-dimensional case (see \cite{Sergeyev et
al. (2001),Sergeyev and Markin (1995),Strongin (1978),Strongin and
Markin (1986),Strongin and Sergeyev (2000)}). Formally, supposing
that both the objective function $f(x)$ and constraints $g_j(x),
1\,\leq j\,\leq m,$ satisfy the Lipschitz condition and the
feasible region is not empty, this problem can be formulated as
follows.

It is necessary to find a point $x^*$ and the corresponding
value\, $g_{m+1}^*$\, such that
 \beq
  g_{m+1}^*=g_{m+1}(x^*)=\min\{g_{m+1}(x)\,:x\in Q_{m+1}\}, \label{6}
 \eeq
where, in order to unify the description process,   the
designation $g_{m+1}(x) \triangleq f(x)$ has been used and regions
$Q_j, 1\,\leq j\,\leq m+1,$ are defined by the rules
 \beq
  Q_1\,=\,[a,b],\,\,\,\,\mbox
  Q_{j+1}\,=\,\{x\,\in\,Q_j\,:\,\textit{g}_j(x)\,\leq
  0\},\,\,\,\,1\,\leq j\,\leq m,  \label{2}
 \eeq
$$
 Q_1\,\supseteq Q_2\,\supseteq \ldots \supseteq
 Q_m\,\supseteq\,Q_{m+1}.
$$
Note that since the constraints $g_j(x)$, $1\leq j\leq m$, are
multiextremal, the admissible region $Q_{m+1}$ and regions
$Q_{j},$ $1\le j\le m,$ can be collections of several disjoint
subregions. We suppose hereafter that all of them consist of
intervals of a finite length.

We assume also that the functions $g_j(x),1\leq j\leq m+1$,
satisfy the corresponding Lipschitz conditions
 \beq
  \mid g_j(x')\,-\,g_j(x'') \mid \, \leq L_j\,\mid
  x'\,-\,x'' \mid,\,\,\,x',\,x''\,\in\,Q_j,\,\,\,1\leq
  j\,\leq m+1, \label{4} \eeq \beq
  0\,<\, L_j\,<\,\infty,\,\,\,\,\,\,1\leq j\leq m+1. \label{5}
 \eeq

\begin{figure}[tp]
\caption{Problem number~6 from \cite{Famularo et al. (2002)} where
functions $f(x)$ and $g_1(x),g_2(x)$ are defined over the whole
search region $[0 , 1.5 \, \pi]$} \label{Fig6_2}
    \epsfig{ figure = 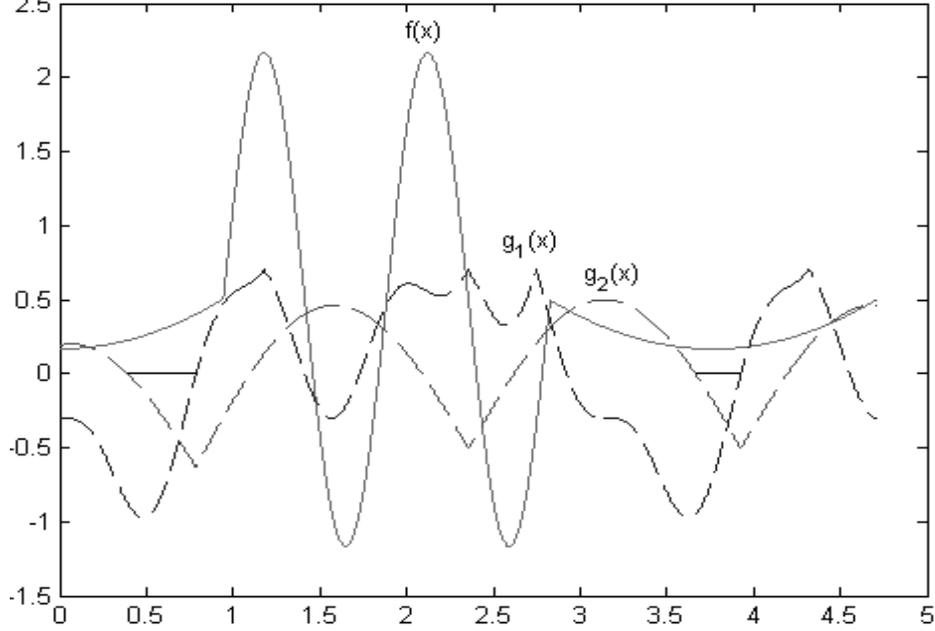, width = 4.8in, height = 3.3in,  silent = yes }
\end{figure}

In order to illustrate the problem under consideration and to
highlight its difference with respect to problems where
constraints and the objective function are defined over the whole
search region, let us consider an example -- test problem number~6
from \cite{Famularo et al. (2002)} shown in Figure~\ref{Fig6_2}.
The problem has two multiextremal constraints and is formulated as
follows \beq
  f^*=f(x^*)=\min\{f(x): g_{1}(x) \le 0, g_{2}(x) \le 0, x\in [0 , 1.5 \, \pi] \},  \label{example_-1}
 \eeq
where
 \beq
  g_{3}(x) \triangleq  f(x)   = \left\{ \begin{array}{ll} \displaystyle \frac{1}{3} \left( \frac{100}{9
    \pi^2} x^2 + \frac{1}{2} \right), & x \le \frac{3 \pi}{10}, \\[12pt]
    \displaystyle  \frac{5}{3} \, \sin\left( \frac{20}{3} x \right)+ \frac{1}{2}, &
    \frac{3 \pi}{10} < x \le \frac{9 \pi}{10}, \\[12pt]
    \displaystyle \frac{1}{3} \left( \frac{100}{9 \pi^2} x^2-\frac{80}{3 \pi}
    x+\frac{33}{2} \right), & x > \frac{9 \pi}{10},
 \end{array} \right.   \label{example_3}
 \eeq
 \beq
  g_1(x) = \displaystyle
 \frac{7}{10}-\left|\sin^3(3 x)+\cos^3(x)\right|,  \label{example_4}
 \eeq
  \beq
  g_2(x) = \displaystyle
 -\left|\frac{(x-\pi)^3}{100}\right|+\left|\cos(2
 (x-\pi))\right|-\frac{1}{2}.  \label{example_5}
 \eeq
The admissible region of problem
(\ref{example_-1})--(\ref{example_5}) consists of two disjoint
subregions shown in Figure~\ref{Fig6_2} at the line $f(x)=0$; the
global minimizer is $x^{*}=3.76984$.

\begin{figure}[tp]
\caption{Graphical representation of problem
(\ref{example_0})--(\ref{example_2})} \label{Fig6_1}
   \epsfig{ figure = 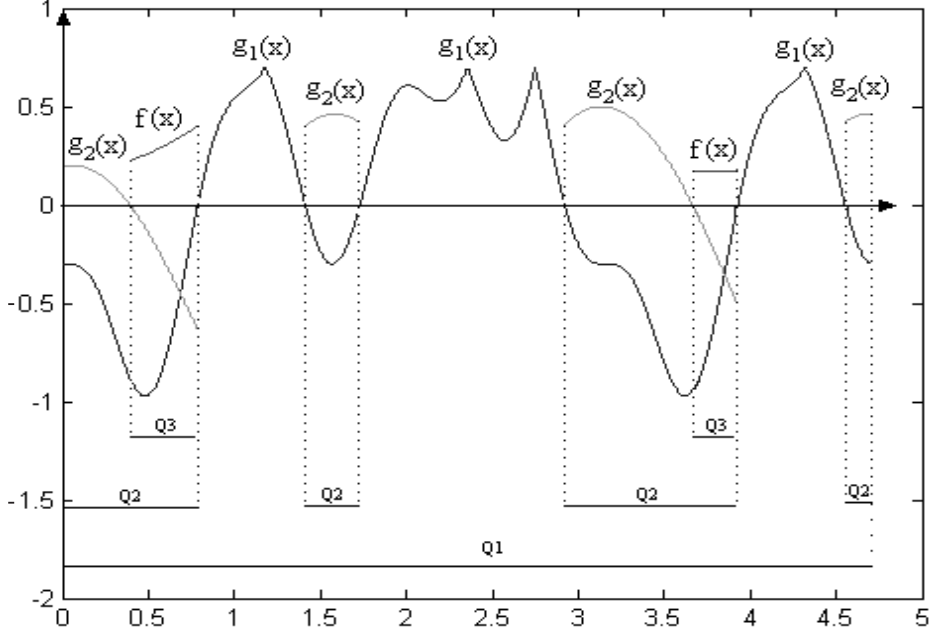, width = 4.8in, height = 3.3in,  silent = yes }
\end{figure}

Problem of the type (\ref{6})--(\ref{5}) considered in this paper
and using the same functions $g_{1}(x)$--$g_{3}(x)$ from
(\ref{example_3})--(\ref{example_5}) is shown in
Figure~\ref{Fig6_1}. It has the same global minimizer
$x^{*}=3.76984$ and is formulated as follows
 \beq
  Q_1\,=\,[0 , 1.5 \, \pi],\,\,\,\,
  Q_{2}\,=\,\{x\,\in\,Q_1 :\,g_1(x)\,\leq
  0\}, \label{example_0}
 \eeq
\beq
 Q_{3}\,=\,\{x\,\in\,Q_2 :\,g_2(x)\,\leq
  0\},  \label{example_1}
 \eeq
\beq
  g_{3}^*=g_{3}(x^*)=\min\{g_{3}(x):x\in Q_{3}\}.  \label{example_2}
 \eeq
It can be seen from Figure~\ref{Fig6_1} that both $g_{2}(x)$   and
  $f(x)$ are partially defined: $g_{2}(x)$ is defined only over
$Q_{2}$ and the objective function $f(x)$ is defined only over
$Q_{3}$ which coincides with the admissible region of problem
(\ref{example_-1})--(\ref{example_5}).

It is not easy to find a traditional algorithm for solving problem
(\ref{6})--(\ref{5}). For example, the penalty approach  requires
that $f(x)$ and $g_i(x),\, 1\le i\le m,$ are defined over the
whole search interval $[a,b]$. At first glance it seems that at
the regions where a function is not defined it can be simply
filled in with either a big number or the function value at the
nearest feasible point. Unfortunately, in the context of Lipschitz
algorithms, incorporating such ideas can lead to infinitely high
Lipschitz constants, causing degeneration of the methods and
non-applicability of the penalty approach.

A promising approach called the {\it index scheme} has been
proposed in \cite{Strongin (1984)} (see also \cite{Sergeyev and
Markin (1995),Strongin and Markin (1986),Strongin and Sergeyev
(2000)}) in combination with information stochastic Bayesian
algorithms for solving problem (\ref{6})--(\ref{5}). An important
advantage of the index scheme is that it does not introduce
additional variables and/or parameters as traditional  approaches
do (see, e.g, \cite{Bertsekas (1996),Horst and Pardalos
(1995),Nocedal and Wright (1999)}). It has been recently shown in
\cite{Sergeyev et al. (2001)} that the index scheme can be also
successfully used in combination with the Branch-and-Bound
approach if the Lipschitz constants $L_j,1\leq j\leq m+1,$ from
(\ref{4}), (\ref{5}) are known a priori.

However, in practical applications (see, e.g. \cite{Pinter
(1996)}) the Lipschitz constants $L_j,1\leq j\leq m+1,$ are very
often unknown. Thus, the problem of their estimating arises
inevitably. If there exists an additional information allowing us
to obtain  a priori fixed constants~$K_{j}$, $1 \le j \le m+1$,
such that
 \[
  L_{j}  < K_{j} < \infty,
  \hspace{3mm} 1\le j\le m+1,  
 \]
then the algorithm IBBA from \cite{Sergeyev et al. (2001)} can be
used.

In this paper,  the case where there is no any additional
information about the Lipschitz constants is considered. A new GO
Algorithm with Local Tuning (ALT) adaptively estimating the local
Lipschitz constants  during the search is proposed. The local
tuning technique introduced in \cite{Sergeyev (1995a),Sergeyev
(1995b),Sergeyev (1998)} for solving unconstrained problems allows
one to accelerate the search significantly in comparison with the
methods using estimates of the global Lipschitz constant. The new
method ALT unifies this approach with  the index scheme  and
geometric ideas allowing one to construct auxiliary functions
similar to minorants used in the IBBA (see \cite{Sergeyev et al.
(2001)}). In a series of numerical experiments it is shown that
usage of adaptive local estimates calculated during the search
instead of a priori given estimates of global Lipschitz constants
accelerates the search significantly.

\section{A New Geometric Index Algorithm with Local Tuning}

Let us associate with every point of the interval $[a,b]$ an
\textit{index}
$$
 \nu\,=\,\nu (x),\,\,\,\, 1\,\leq \nu\,\leq M,
$$
which is defined by the conditions
 \beq
   g_j(x)\,\leq 0,\,\,\,\, 1\,\leq j\,\leq\,\nu\,-\,1,\,\,\,\,\,
   g_{\nu}(x)\,>\,0,  \label{7}
 \eeq
where for $\nu\,=\,m+1$ the last inequality is omitted. We shall
call a \textit{trial} the operation of evaluation of the functions
$g_j(x),\,1\,\leq\,j\,\leq\,\nu(x)$, at a point $x$.

Thus, the index scheme considers constraints one at a time at
every point where it has been decided to try to calculate the
objective function $g_{m+1}(x)$. Each constraint~$g_i(x)$ is
evaluated at a point $x$ only if all the inequalities
 $$
  g_j(x)\,\leq 0,\,\,\,\,1\,\leq j\,< i,
 $$
have been satisfied at this point. In its turn the objective
function\, $g_{m+1}(x)$ is computed only for those points where
all the constraints have been satisfied.

Suppose now that $k+1$, $k \geq 1$, trials have been executed at
some points
 \beq
  a=x_0 < x_1 < \ldots < x_i < \ldots < x_k =b \label{10}
 \eeq
and the index\, $\nu_i=\nu(x_i),\,\,0\leq i \leq k,$\, have been
calculated following $(\ref{7})$. Due to the index scheme, the
estimate
 \beq
  z^*_k =\displaystyle\ \min \,\{ g_{M^k}(x_i):\, 0 \leq i \leq k ,\, \nu(x_i) =
  M^k\}   \label{12}
 \eeq
of the minimal value of the function $g_{M^k}(x)$ found after $k$
iterations can be calculated and the values
 \beq
 z_i=   g_{\nu(x_i)}(x_i) - \left\{
 \begin{array}{cl}
 0 ,& \hspace{2mm} \mbox{if} \,\,\, \nu(x_i)< M^k\\
 z^*_k ,& \hspace{2mm} \mbox{if} \,\,\, \nu(x_i)=M^k\\
\end{array}\right.   \label{11}
 \eeq
can be associated with the points $x_i$ from~(\ref{10}).

In the new algorithm ALT, we propose to adaptively estimate at
each iteration   the local Lipschitz constants over subintervals
$[x_{i-1}, x_i] \subset [a,b]$, $1 \leq i \leq k$, by using the
information obtained from executing trials at the points $x_{i}$,
$0 \leq i \leq k$, from (\ref{10}). Particularly, at each point
$x_i$, $0 \leq i \leq k,$ having the index~$\nu_{i}$ we calculate
a local estimate $\eta_{i}$ of the Lipschitz constant
$L_{\nu_{i}}$ at a neighborhood of the point $x_i$ as follows
 \beq
   \eta_{i} =  \max \{ \lambda_i, \gamma_i, \xi\}, \hspace{5mm} 0 \le i \le
   k.
   \label{eta}
 \eeq
Here, $\xi > 0$ is a small number reflecting our supposition that
the objective function and constraints are not just constants over
$[a,b]$, i.e., $L_{j} \ge \xi$, $1 \le j \le m+1$. The values
$\lambda_{i}$ are calculated as follows
 \beq
  \hspace{-3mm}\small{ \lambda_i   =
 \left\{\begin{array}{ll}
 \max\{\mid z_j - z_{j-1} \mid (x_j - x_{j-1})^{-1}:\,\, j=i, i+1\}, & \mbox{if}\,\,\nu_{i-1}=\nu_i=\nu_{i+1} \\
 & \\
 \max\{\mid z_i - z_{i-1} \mid (x_i - x_{i-1})^{-1},\,\,\,z_i(x_{i+1} - x_i)^{-1} \}, & \mbox{if}\,\,\nu_{i-1}=\nu_i<\nu_{i+1}\\
 & \\
  \max\{\mid z_{i+1} - z_i \mid (x_{i+1} - x_i)^{-1},\,\,\,z_i(x_i - x_{i-1})^{-1}\}, & \mbox{if}\,\,\nu_{i-1}>\nu _i=\nu_{i+1} \\
 & \\
  \max\{z_i (x_i - x_{i-1})^{-1},\,\,\,z_i (x_{i+1} - x_i)^{-1}\}, & \hspace{-5.5mm} \mbox{if}\,\,\nu_i< \nu_{i-1},\,\,\nu_i<\nu_{i+1}\\
 & \\ \label{lambda}
  z_i (x_i - x_{i-1})^{-1}, &  \mbox{if}\,\, \nu_{i-1} > \nu_i > \nu_{i+1}\\
 & \\
  z_i (x_{i+1} - x_{i})^{-1}, & \mbox{if}\,\, \nu_{i-1} < \nu_i < \nu_{i+1}\\
 & \\
 \mid z_i - z_{i-1} \mid (x_i - x_{i-1})^{-1}, & \mbox{if}\,\,\nu_{i-1}=\nu_i>\nu_{i+1}\\
 & \\
 \mid z_{i+1} - z_i \mid (x_{i+1} - x_i)^{-1}, & \mbox{if}\,\,\nu_{i-1}<\nu_i=\nu_{i+1}\\
 & \\
  0,        & \mbox{otherwise}                      \\
 \end{array} \right. }
\eeq where $z_i$, $0 \leq i \leq k$, are from (\ref{11}).
Naturally, when $i = 0$ or $i = k$, only one of the two
expressions in the first four cases are defined and are used to
calculate $\lambda_i$. The values $\gamma_i$, $0 \leq i \leq k,$
are calculated in the following way:
 \beq
   \gamma_i = \Lambda_{\nu_i} \max\{ x_i - x_{i-1}, x_{i+1} -
   x_{i}\}/X_{\nu_i}^{\max},
   \label{gamma}
 \eeq
 \beq
  \Lambda_{\nu_i} =  \Lambda_{\nu_i} (k) =
  \max\{ \Lambda_{\nu_i}(k-1), \max \{\lambda_j: \nu_j = \nu_i,\,0\leq j \leq
  k\}\}, \label{Lambda}
 \eeq
where $\Lambda_{\nu_i}$ are adaptive estimates of the global
Lipschitz constants $L_{\nu_i}$ and
 \beq
   X_{\nu_i}^{\max}\,=\,\max\{x_j\,-\,x_{j-1}\,:\,\nu_j=\nu_i\,\,
   {\rm or} \,\,\nu_{j-1}=\nu_i,\,\,1\leq j \leq k\}.
 \eeq

The values $\lambda _{i}$ and $\gamma _{i}$ reflect the influence
on $\eta_i$ of the local and global information obtained during
the previous iterations. When both intervals $[x_{i-1},x_{i}]$ and
$[x_{i},x_{i+1}]$ are small, then $\gamma _{i}$ is small too (see
(\ref{gamma})) and, due to (\ref{eta}), the local information
represented by $\lambda _{i}$ has  major importance. The value
$\lambda _{i}$ is calculated by considering the intervals
$[x_{i-2},x_{i-1}], \; [x_{i-1},x_{i}],$ and $ [x_{i},x_{i+1}]$
(see (\ref{lambda})) as those which have the strongest influence
on the local estimate at the point $x_i$ and, in general, at the
interval $[x_{i-1},x_{i}]$. When at least one of the intervals
$[x_{i-1},x_{i}]$, $[x_{i},x_{i+1}]$ is very wide, the local
information is not reliable and, due to (\ref{eta}), the global
information represented by $\gamma _{i}$ has the major influence
on $\eta_{i}$. Thus, local and global information are balanced in
the values $\eta_{i}$, $0\le i\le k$. Note that the method uses
the local information over the {\it whole} search region $[a,b]$
{\it during} the global search both for the objective function and
constraints.

We are ready now to describe the new algorithm ALT.

\begin{description}
\item[\textbf{Step 0 (Initialization).}] Suppose that $k+1$, $k \ge 1$, trails have been
already executed in a way at points
 \beq
  x^{0}=a, x^{1}=b, x^{2}, x^{3}, ... ,x^{i}, ...
  x^{k-1}, x^{k} \label{17}
 \eeq
and their indexes and the value
 \beq
   M^{k}=\max\{\nu(x^i):0\leq\,i\,\leq\,k\}  \label{18}
 \eeq
have been calculated. The value $M^{k}$  defined in (\ref{18}) is
the maximal index  obtained  during the search after $k+1$ trials.
The choice of the point $x^{k+1}$, $k\geq\,1$, where the next
trial will be executed is determined by the rules presented below.

\item[\textbf{Step 1.}]
Renumber the points $x^{0},....,x^{k}$ of the previous $k$
iterations by subscripts\footnote{Thus, two numerations are used
during the work of the algorithm. The record $x^{i}$ from
(\ref{17}) means that this point has been generated during the
$i$-th iteration of the ALT. The record $x_i$ indicates the place
of the point in the row (\ref{10}). Of course, the second
enumeration is changed during every iteration.} in order to form
the sequence (\ref{10}).

\item[\textbf{Step 2.}]
Recalculate the estimate $z^*_k$ of the minimal value of the
function $g_{M^k}(x)$ found after $k$ iterations and the values
$z_i$ by using formulae (\ref{12}) and (\ref{11}), respectively.
For each trial point $x_i$ having the index $\nu_i$, $0\leq i \leq
k$, calculate estimate $\eta_i$  from (\ref{eta}).

\item[\textbf{Step 3.}]
For each interval $[x_{i-1}, x_i],\, 1\leq i\leq k$, calculate the
\textit{characteristic} of the interval
 \beq
 \hspace{-4mm}  R_i  =\left\{
         \begin{array}{ll}
 ( \frac{1}{\eta_i+\eta_{i-1}}) \left[\eta_i z_{i-1}+\eta_{i-1}z_i+r\eta_{i-1}\eta_i(x_{i-1}-x_i)\right],& \, \nu_{i-1}\,=\,\nu_i \\
 z_i - r\eta_i(x_i - x_{i-1}\,-z_{i-1}/r\eta_{i-1}), & \, \nu_{i-1}<\,\nu_i\\
 z_{i-1} - r\eta_{i-1}(x_i - x_{i-1}\,-z_i/r\eta_i), & \, \nu_{i-1}>\,\nu_i\\
         \end{array} \right.     \label{19}
 \eeq
where $r>1$ is the reliability parameter of the method (this kind
of parameters is quite traditional in   Lipschitz global
optimization; discussions related to its choice and meaning can be
found in \cite{Pinter (1996),Strongin (1978),Strongin and Sergeyev
(2000)}).

\item[\textbf{Step 4.}] Find an interval $t$ corresponding to the minimal
characteristic, i.e.,
 \beq
   t= \arg \min\{\,R_i:\,1\leq i \leq k \}.   \label{20}
 \eeq
If the minimal value of the characteristic is attained for several
subintervals, then the minimal integer satisfying~(\ref{20}) is
accepted as $t$.

\item[\textbf{Step 5.}]  If for the interval $[x_{t-1}, x_t],$ where $t$ is from (\ref{20}), the stoping rule
 \beq
    x_t-x_{t-1} \le \varepsilon (b-a),  \label{21}
 \eeq
where $a$ and $b$ are from (\ref{6})--(\ref{5}), is satisfied for
a preset accuracy $\varepsilon>0$, then \textbf{Stop} -- the
required accuracy has been reached. In the opposite case, go to
Step~6.

\item[\textbf{Step 6.}] Execute the $(k+1)$-th trial at the point
 \beq
\hspace{-3mm}   x^{k+1}  =\left\{
             \begin{array}{ll}
 (\frac{1}{r\eta_t+r\eta_{t-1}}) \left[z_{t-1}-z_t+r\eta_{t-1}x_{t-1}+r\eta_tx_t\right], \,& {\rm if} \,\nu_{i-1}=\,\nu_i \\
 0.5(x_{t-1} + x_t),& {\rm if} \, \nu_{i-1} \neq \nu_i,
             \end{array} \right. \label{22}
 \eeq
and evaluate its index $\nu(x^{k+1})$.

\item[\textbf{Step 7.}] This step consists of the following
alternatives:

\textit{Case 1.} If $\nu(x^{k+1}) > M^{k}$, then perform two
additional trials at the points
 \beq
   x^{k+2} = 0.5(x_{t-1} + x^{k+1}),   \hspace{1cm}
   x^{k+3} = 0.5(x^{k+1} + x_t),  \label{24}
 \eeq
calculate their indexes, set $k=k+3$, and go to Step 8.

\textit{Case 2.} If $\nu(x^{k+1}) < M^{k}$ and among the points
(\ref{10}) there exists only one point $x_T$ with the maximal
index $M^{k}$, i.e., $\nu_T = M^{k}$, then execute two additional
trials at the points
 \beq
   x^{k+2} = 0.5(x_{T-1} + x_T), \label{25}
 \eeq
 \beq
   x^{k+3} = 0.5(x_T + x_{T+1}), \label{26}
 \eeq
if $0 < T < k$, calculate their indexes, set $k=k+3$,  and go to
Step~8. If $T=0$ then the trial is executed only at the point
(\ref{26}). Analogously, if $T=k$ then the trial is executed only
at the point (\ref{25}). In these two cases, calculate the index
of the additional point, set $k=k+2$ and go to Step 8.

\textit{Case 3.} In all the remaining cases set $k=k+1$ and go to
Step 8.

\item[\textbf{Step 8.}] Calculate $M^{k}$ and go to Step~1.
\end{description}

Global convergence conditions of the ALT are described by the
following two theorems given, due to the lack of space, without
proofs that can be derived using Theorems~2 and~3 from
\cite{Sergeyev (1995b)} and Theorem~2 from \cite{Sergeyev et al.
(2001)}.
\begin{theorem}
Let the feasible region $Q_{m+1} \neq \varnothing$ consists of
intervals having finite lengths, $x^{*}$ be any solution to
problem (\ref{6})--(\ref{5}), and $j=j(k)$ be the number of an
interval $[x_{j-1},x_{j}]$ containing this point during the $k$-th
iteration. Then, if for $k\ge k^{*}$ the following conditions
 \beq
r\Lambda_{\nu_{j-1}}> C_{j-1}, \hspace{5mm}  r\Lambda_{\nu_{j}}>
C_{j}, \label{29}
 \eeq
 \beq
C_{j-1}=         z_{j-1}/(x^{*}-x_{j-1}), \hspace{5mm} C_{j}=
z_{j}/(x_{j}-x^{*}).                \label{30}
 \eeq
take  place, then the point $x^{*}$ will be a limit point of the
trial sequence $\{x^{k}\}$  generated by the ALT.
  \label{t1}
\end{theorem}
\begin{theorem}
For any problem   (\ref{6})--(\ref{5}) there exists a value $r^*$
such that conditions (\ref{29}) are satisfied  for all parameters
$r>r^*,$ where $r$ is from (\ref{19}) and (\ref{22}).
  \label{t2}
\end{theorem}

\section{Numerical Comparison}
The new algorithm has been numerically compared with the following
methods:
\begin{description}
\item[--] The method proposed by Pijavskii (see \cite{Horst and Pardalos (1995),Pijavskii (1972)}) combined with the penalty approach used to
reduce the constrained problem to an unconstrained one; this
method is indicated hereafter as PEN. The Lipschitz constant of
the obtained unconstrained problem is supposed to be known as it
is required by Pijavskii algorithm.

\item[--] The method IBBA from \cite{Sergeyev et al. (2001)}
using the index scheme in combination with the Branch-and-Bound
approach and the known Lipschitz constants $L_j,1\leq j\leq m+1,$
from (\ref{4}), (\ref{5}).
\end{description}

Ten non-differentiable test problems introduced in \cite{Famularo
et al. (2002)} have been used in the experiments (since there were
several misprints in the original paper \cite{Famularo et al.
(2002)}, the accurately verified formulae have been applied, which
are available at the Web-site {\small
http://wwwinfo.deis.unical.it/$\sim$yaro/constraints.html}). In
this set of tests, problems 1--3 have one constraint, problems
4--7 two constraints, and problems 8--10 three constrains.

In these test problems,  all constrains and the objective function
are defined over the whole region $[a,b]$ from (\ref{2}). These
test problems were used because the PEN needs this additional
information for its work and is not able to solve problem
(\ref{6})--(\ref{5}). Naturally, the methods IBBA and ALT solved
all the problems using the statement (\ref{6})--(\ref{5}) and did
not take benefits from the additional information given  (see
examples  from Figures~\ref{Fig6_2} and~\ref{Fig6_1}) by the
statement (\ref{example_-1})--(\ref{example_5}) in comparison with
(\ref{example_0})--(\ref{example_2}).

In order to demonstrate the influence of changing the search
accuracy~$\varepsilon$ on the convergence speed of the methods,
two different values of $\varepsilon$, namely, $\varepsilon =
10^{-4}$ and $\varepsilon = 10^{-5}$ have been used. The same
value $\xi = 10^{-6}$ from~(\ref{eta}) has been used in all the
experiments for all the methods.

Table \ref{PEN} represents the results for the PEN (see
\cite{Famularo et al. (2002)}). The constrained problems were
reduced to the unconstrained ones as follows
 \beq
  f_{P^*}(x)=f(x) + P^* \max\left\{ 0, g_1(x), g_2(x), \dots,
  g_{N_v}(x)\right\}.  \label{pen}
 \eeq
The column ``Eval.'' in Table~\ref{PEN} shows the total number of
evaluations of the objective function $f(x)$ and all the
constraints. Thus, it is equal to
 \[  (N_v +1)\times N_{trials}, \]
where $N_v$ is the number of constraints and $N_{trials}$ is the
number of the trials executed by the PEN for each problem.

Results obtained by the IBBA (see \cite{Sergeyev et al. (2001)})
and by the new method ALT with the parameter $r=1.3$ are
summarized in Tables~\ref{IBBA} and ~\ref{ALT_1.3}, respectively.
Columns in the tables have the following meaning for each value of
the search accuracy $\varepsilon$:

$-$the column $\mathcal{N}$ indicates the problem number;

$-$the columns $N_{g_1}, N_{g_2},$ and $N_{g_3}$ represent the
number of trials where the constraint $g_i,\, 1\leq i\leq 3$, was
the last evaluated constraint;

$-$the column ``Trials'' is the total number of trial points
generated by the methods;

$-$the column ``Eval.'' is the total number of evaluations of the
objective function and the constraints. This quantity is equal to:

$-$$N_{g_1} + 2 \times  N_f$, for problems with one constraint;

$-$$N_{g_1} + 2 \times  N_{g_2} + 3 \times N_f$, for problems with
two constraints;

$-$$N_{g_1} + 2 \times  N_{g_2} + 3 \times N_{g_3} + 4 \times
N_f$, for problems with three constraints.

\begin{table}[t]
 \begin{center}
\caption{Numerical results obtained by the PEN}
    \label{PEN}
    \begin{tabular}{@{\extracolsep{\fill}}c c | r r | r r}
\hline \raisebox{-1.5ex}[0ex][0ex]{$\mathcal{N}$} &
       \raisebox{-1.5ex}[0ex][0ex]{$P^*$}
&\multicolumn{2}{c |}{$\varepsilon = 10^{-4}$} &\multicolumn{2}{c}{$\varepsilon = 10^{-5}$} \\
& & Trials & Eval. & Trials & Eval. \\
\noalign{\smallskip}\hline\noalign{\smallskip}
  1  & $15$ & $247$ & $494$ & $419$ & $838$ \\
  2  & $15$ & $241$ & $482$ & $313$ & $626$ \\
  3  & $15$ & $917$ & $1834$ & $2127$ & $4254$ \\
  4  & $15$ & $273$ & $819$  & $861$ & $2583$ \\
  5  & $20$ & $671$ & $2013$  & $1097$ & $3291$ \\
  6  & $15$ & $909$ & $2727$ & $6367$ & $19101$ \\
  7  & $15$ & $199$ & $597$   & $221$ & $663$ \\
  8  & $15$ & $365$ & $1460$  & $415$ & $1660$ \\
  9  & $15$ & $1183$ & $4732$ & $4549$ & $18196$ \\
  10 & $15$ & $135$ & $540$   & $169$ & $676$ \\
  Av. & $-$ & $514.0$ & $1569.8$ & $1653.8$ & $5188.8$ \\
  \noalign{\smallskip}\hline
    \end{tabular}
  \end{center}
\end{table}

The asterisk in Table~\ref{ALT_1.3} indicates that $r=1.3$ was not
sufficient to find the global minimizer of problem~7. The results
for this problem in Table~\ref{ALT_1.3} are obtained using the
value $r=1.9$; the ALT with this value finds the solution.
Finally, Table~\ref{speed_1.3} represents the improvement (in
terms of the number of trials and evaluations) obtained by the ALT
in comparison with the other methods used in the experiments.

As it can be seen from Tables~\ref{PEN}--\ref{speed_1.3}, the
algorithms IBBA and ALT constructed in the framework of the index
scheme significantly outperform the traditional method PEN. The
ALT demonstrates a high improvement in terms of the trials
performed with respect to the IBBA as well. In particular, the
greater the difference between estimates of the local Lipschitz
constants (for the objective function or for the constraints), the
higher is the speed up obtained by the ALT (see
Table~\ref{speed_1.3}). The improvement is especially high if the
global minimizer lies inside of a feasible subregion with a small
(with respect to the global Lipschitz constant) value  of the
local Lipschitz constant as it happens for example for problems 6
and 9. The advantage of the new method is more pronounced when the
search accuracy $\varepsilon$ increases (see
Table~\ref{speed_1.3}).

\begin{table}[tp]
 \begin{center}
 \caption{Numerical results obtained by the IBBA}
    \label{IBBA}
\begin{tabular}{@{\extracolsep{-0.25mm}}c | r r r r r r | r r r r r r}
\hline \raisebox{-1.5ex}[0ex][0ex]{$\mathcal{N}$}
 &\multicolumn{6}{c |}{$\varepsilon = 10^{-4}$} &\multicolumn{6}{c}{$\varepsilon = 10^{-5}$} \\
 &$N_{g_1}$\hspace{-1mm} & $N_{g_2}$\hspace{-1mm} &
$N_{g_3}$\hspace{-1mm} & $N_{f}$ & Trials & Eval.
 &$N_{g_1}$\hspace{-1mm} & $N_{g_2}$\hspace{-1mm} &
$N_{g_3}$\hspace{-1mm} & $N_{f}$ & Trials & Eval. \\
\noalign{\smallskip}\hline\noalign{\smallskip}
 1 & $23$  & $-$  & $-$ & $28$ & $51$  & $79$  & $23$   & $-$   & $-$  & $34$  & $57$  & $91$ \\
 2 & $18$  & $-$  & $-$ & $16$ & $34$  & $50$  & $20$   & $-$   & $-$  & $22$  & $42$  & $64$ \\
 3 & $171$ & $-$  & $-$ & $19$ & $190$ & $209$ & $175$  & $-$   & $-$  & $21$  & $196$ & $217$ \\
 4 & $136$ & $15$ & $-$ & $84$ & $235$ & $418$ & $170$  & $15$  & $-$  & $226$ & $411$ & $878$ \\
 5 & $168$ & $91$ & $-$ & $24$ & $283$ & $422$ & $188$  & $101$ & $-$  & $26$  & $315$ & $468$ \\
 6 & $16$  & $16$ & $-$ & $597$& $629$ & $1839$& $17$   & $17$  & $-$  & $2685$& $2719$& $8106$ \\
 7 & $63$  & $18$ & $-$ & $39$ & $120$ & $216$ & $65$   & $19$  & $-$  & $43$  & $127$ & $232$ \\
 8 & $29$  & $11$ & $3$ & $21$ &  $64$ & $144$ & $29$   & $14$  & $3$  & $23$  &  $69$ & $158$ \\
 9 & $8$   & $86$ & $57$& $183$& $334$ & $1083$& $10$   & $88$  & $57$ & $851$ & $1006$& $3761$ \\
10 & $42$  & $3$  & $17$& $13$ &  $75$ & $151$ & $42$   & $3$   & $17$ & $15$  &  $77$ & $159$\\
Av.& $-$ & $-$  & $-$ & $-$  & $201.5$ & $461.1$ & $-$ & $-$  & $-$ & $-$  & $501.9$ & $1413.4$\\
 \noalign{\smallskip}\hline
 \end{tabular}
 \end{center}


 \begin{center}
 \caption{Numerical results obtained by the ALT with $r=1.3$}
    \label{ALT_1.3}
\begin{tabular}{c | r r r r r r | r r r r r r}
\hline \raisebox{-1.5ex}[0ex][0ex]{$\mathcal{N}$}
 &\multicolumn{6}{c |}{$\varepsilon = 10^{-4}$} &\multicolumn{6}{c}{$\varepsilon = 10^{-5}$} \\
 &$N_{g_1}$\hspace{-1mm} & $N_{g_2}$\hspace{-1mm} &
$N_{g_3}$\hspace{-1mm} & $N_{f}$ & Trials & Eval.
 &$N_{g_1}$\hspace{-1mm} & $N_{g_2}$\hspace{-1mm} &
$N_{g_3}$\hspace{-1mm} & $N_{f}$ & Trials & Eval. \\
\noalign{\smallskip}\hline\noalign{\smallskip}
1   & $27$ & $-$   & $-$  & $17$ & $44$   & $61$   & $27$ & $-$  & $-$  & $19$ & $46$  & $65$ \\
2   & $19$ & $-$   & $-$  & $15$ & $34$   & $49$   & $22$ & $-$  & $-$  & $16$ & $38$  & $54$ \\
3   & $12$ & $-$   & $-$  & $9$  & $21$   & $30$   & $14$ & $-$  & $-$  & $10$ & $24$  & $34$ \\
4   & $45$ & $11$  & $-$  & $37$ & $93$   & $178$  & $45$ & $11$ & $-$  & $48$ & $104$ & $211$ \\
5   & $73$ & $44$  & $-$  & $15$ & $132$  & $206$  & $76$ & $44$ & $-$  & $17$ & $137$ & $215$ \\
6   & $21$ & $11$  & $-$  & $42$ & $74$   & $169$  & $21$ & $11$ & $-$  & $64$ & $96$  & $235$ \\
7$^*$ & $34$ & $27$  & $-$  & $39$ & $100$  & $205$ & $34$ & $34$ & $-$ & $42$ & $110$ & $228$ \\
8   & $12$ & $20$  & $4$  & $23$ & $59$   & $156$  & $12$ & $22$ & $4$  & $24$ & $62$  & $164$ \\
9   & $8$  & $16$  & $3$  & $29$ & $56$   & $165$  & $8$  & $16$ & $3$  & $36$ & $63$  & $193$ \\
10  & $14$ & $2$   & $13$ & $13$ & $42$   & $109$  & $14$ & $2$  & $13$ & $18$ & $47$  & $129$ \\
Av. & $-$  & $-$   & $-$  & $-$  & $65.5$ & $132.8$ & $-$ & $-$  & $-$  & $-$  & $72.7$& $152.8$ \\
 \noalign{\smallskip}\hline
 \end{tabular}
 \end{center}


 \begin{center}
 \caption{Improvement obtained by the ALT with $r = 1.3$ in
comparison with the other methods used in the experiments}
 \label{speed_1.3}
\begin{tabular}{@{\extracolsep{1.5mm}}c | r r r r | r r r r} \hline\noalign{\smallskip}
\raisebox{-1.5ex}[0ex][0ex]{$\mathcal{N}$}
&\multicolumn{4}{c |}{$\varepsilon=10^{-4}$} & \multicolumn{4}{c}{$\varepsilon=10^{-5}$} \\
 &\multicolumn{2}{c}{Trials} & \multicolumn{2}{c |}{Eval.}
 &\multicolumn{2}{c}{Trials} & \multicolumn{2}{c}{Eval.}\\
&{$\frac{\mbox{PEN}}{\mbox{ALT}}$}
&{$\frac{\mbox{IBBA}}{\mbox{ALT}}$}\hspace{-1mm}
&{$\frac{\mbox{PEN}}{\mbox{ALT}}$}
&{$\frac{\mbox{IBBA}}{\mbox{ALT}}$}\hspace{-1mm}
&{$\frac{\mbox{PEN}}{\mbox{ALT}}$}
&{$\frac{\mbox{IBBA}}{\mbox{ALT}}$}\hspace{-1mm}
&{$\frac{\mbox{PEN}}{\mbox{ALT}}$}
&{$\frac{\mbox{IBBA}}{\mbox{ALT}}$} \\
\noalign{\smallskip}\hline\noalign{\smallskip}
1       & $5.61 $ & $1.16 $ & $8.10 $ & $1.30 $ & $9.11 $ & $1.24 $ & $12.89 $ & $1.40 $ \\
2       & $7.09 $ & $1.00 $ & $9.84 $ & $1.02 $ & $8.24 $ & $1.11 $ & $11.59 $ & $1.19 $ \\
3       & $43.67 $ & $9.05 $ & $61.13 $ & $6.97 $ & $88.63 $ & $8.17 $ & $125.12 $ & $6.38 $ \\
4       & $2.94 $ & $2.53 $ & $4.60 $ & $2.35 $ & $8.28 $ & $3.95 $ & $12.24 $ & $4.16 $ \\
5       & $5.08 $ & $2.14 $ & $9.77 $ & $2.05 $ & $8.01 $ & $2.30 $ & $15.31 $ & $2.18 $ \\
6       & $12.28 $ & $8.50 $ & $16.14 $ & $10.88 $ & $66.32 $ & $28.32 $ & $81.28$ & $34.49 $ \\
7$^*$\hspace{-1.5mm }
        & $1.99 $ & $1.20 $ & $2.91 $ & $1.05 $ & $2.01 $ & $1.15 $ & $2.91 $ & $1.02 $ \\
8       & $6.19 $ & $1.08 $ & $9.36 $ & $0.92 $ & $6.69 $ & $1.11 $ & $10.12 $ & $0.96 $ \\
9       & $21.13 $ & $5.96 $ & $28.68 $ & $6.56 $ & $72.21 $ & $15.97 $ & $94.28 $ & $19.49 $ \\
10      & $3.21 $ & $1.79 $ & $4.95 $ & $1.39 $ & $3.60 $ & $1.64 $ & $5.24 $ & $1.23 $ \\
Av.     & $10.92 $ & $3.44 $ & $15.55 $ & $3.45 $ & $27.31 $ & $6.50 $ & $37.10 $ & $7.25 $ \\

 \noalign{\smallskip}\hline
 \end{tabular}
 \end{center}
 \end{table}

\begin{figure}[tp]
\caption{Solving  by the method PEN the unconstrained problem
(\ref{pen}) constructed from the problem
(\ref{example_-1})--(\ref{example_5}) shown in
Figure~\ref{Fig6_2}} \label{f1}
   \epsfig{ figure = 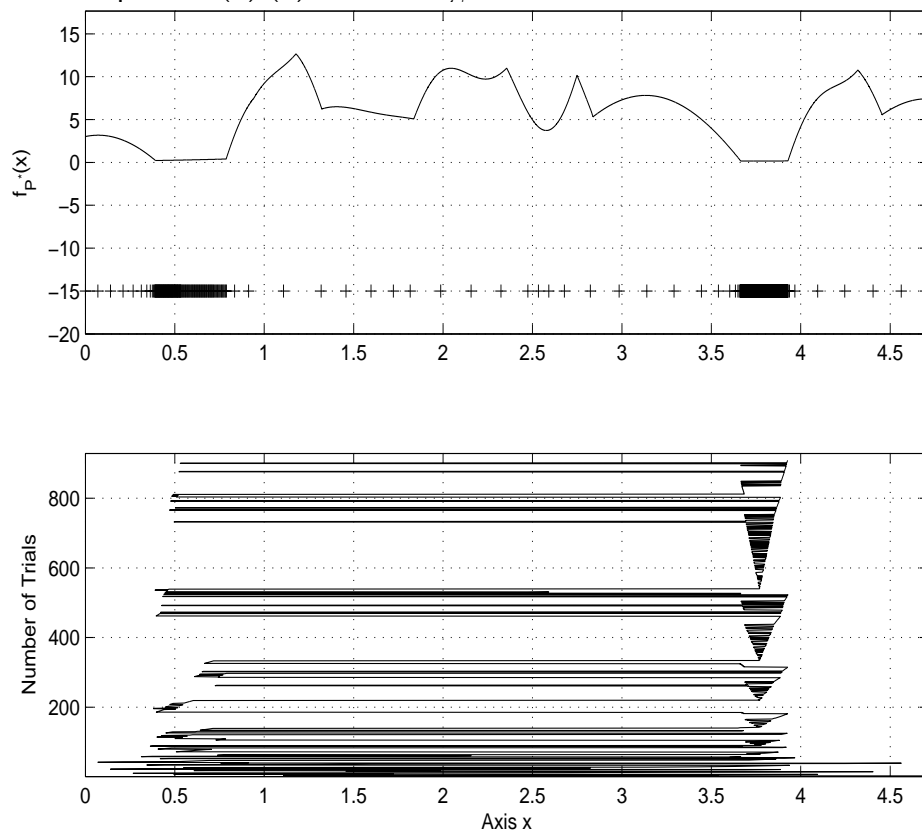, width = 4.8in, height = 4.3in,  silent = yes }
\end{figure}

\begin{figure}[tp]
\caption{Solving the problem (\ref{example_0})--(\ref{example_2})
by the method IBBA} \label{f2}
  \epsfig{ figure = 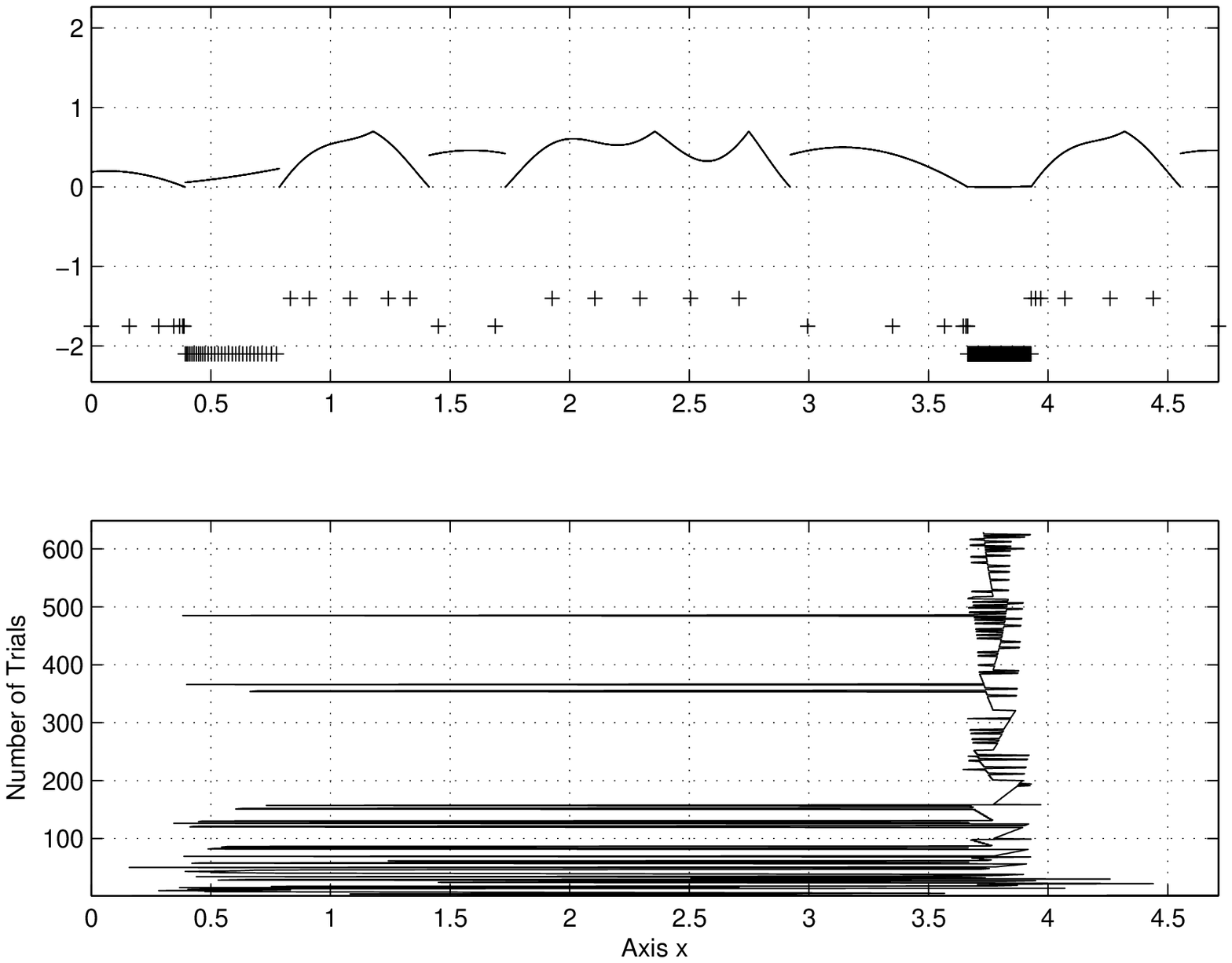, width = 4.8in, height = 4.3in,  silent = yes }

\vspace{5mm}

\caption{Solving the problem (\ref{example_0})--(\ref{example_2})
by the method ALT with $r=1.3$} \label{f3}
   \epsfig{ figure = 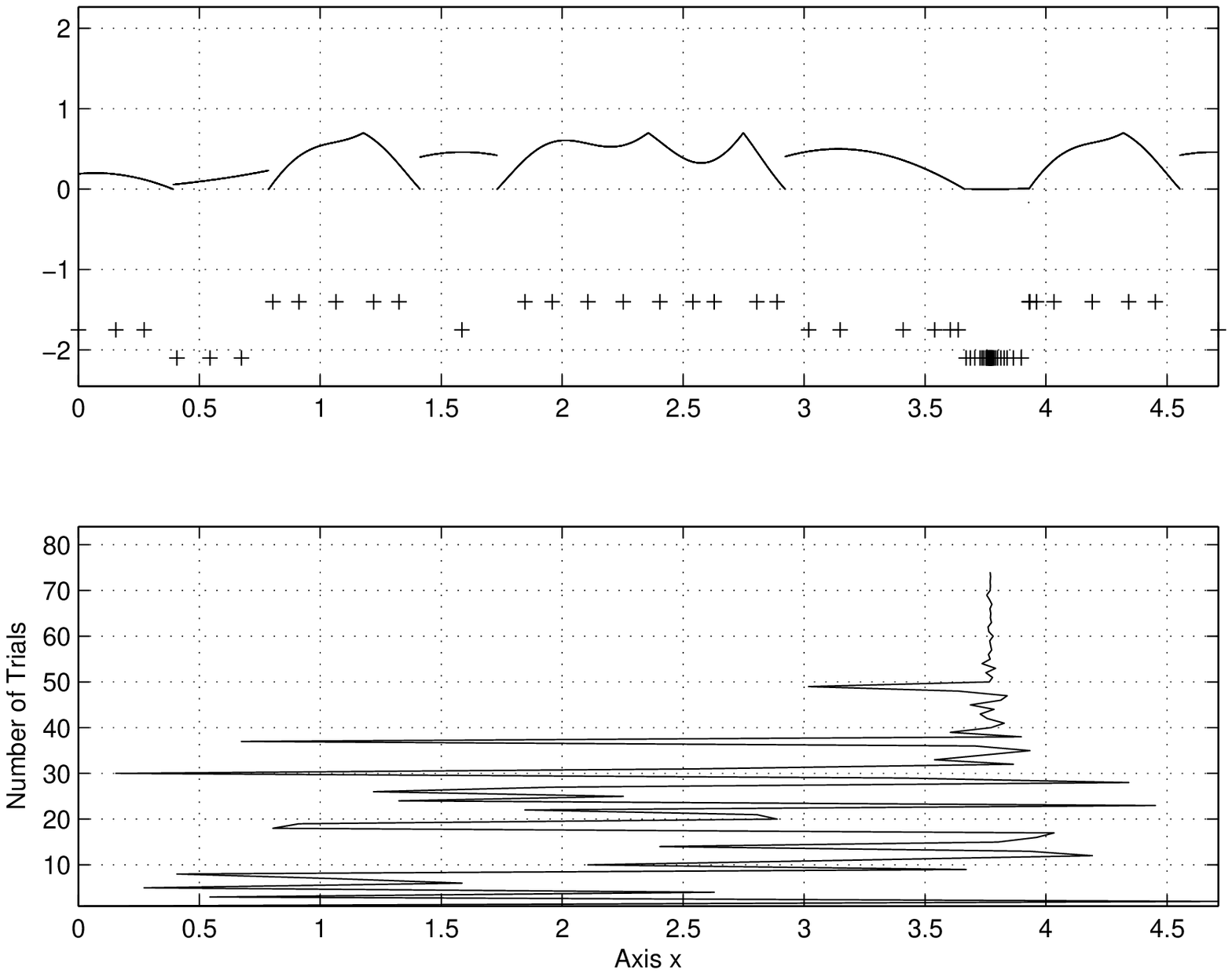, width = 4.8in, height = 4.3in,  silent = yes }
\end{figure}

In order to illustrate performance of the methods graphically, in
Figures~\ref{f1}--\ref{f3} we show dynamic diagrams of the search
(with the accuracy $\varepsilon = 10^{-4}$ in~(\ref{21})) executed
by the PEN, the IBBA, and the ALT, respectively, for problem~6
from \cite{Famularo et al. (2002)}. The upper subplot of
Figure~\ref{f1} contains the function $f_{P^*}(x)$ from
(\ref{pen}) constructed from the problem
(\ref{example_-1})--(\ref{example_5}) shown in
Figure~\ref{Fig6_2}. The upper subplots of Figures~\ref{f2} and
\ref{f3} contain the index function (see \cite{Sergeyev et al.
(2001),Strongin and Sergeyev (2000)} for a detailed discussion)
corresponding to the problem (\ref{example_0})--(\ref{example_2})
from Figure~\ref{Fig6_1}. Note that the local Lipschitz constant
corresponding to the objective function over this subregion is
significantly smaller than the global one (see
Figures~\ref{Fig6_1} and~\ref{f3}).

The line of symbols `+' located under the graph of the function
(\ref{pen}) in   Figure~\ref{f1} shows  points at which trials
have been executed by the PEN. The lower subplots show dynamics of
the search. The PEN has executed 909 trials and the number of
evaluations was equal to $909 \times 3 = 2727$. In Figure~\ref{f2}
the first line (from up to down) of symbols `+', located under the
graph of problem~(\ref{example_0})--(\ref{example_2}), represents
the points where the first constraint has not been satisfied
(number of such trials is equal to~16). Thus, due to the decision
rule of the IBBA, the second constraint has not been evaluated at
these points.

The second line of symbols `+' represents the points where the
first constraint has been satisfied but the second constraint has
been not (number of such trials is equal to 16). At these points
both constraints have been evaluated but the objective function
has been not. The last line represents the points where both
constraints have been  satisfied (number of such trials is 597)
and, therefore, the objective function has been evaluated too. The
total number of evaluations is equal to $16+16 \times 2+597 \times
3 = 1839$. These evaluations have been executed during $16 + 16 +
597 = 629$ trials.

Similarly, in Figure~\ref{f3}, the first line of symbols `+'
indicates 21 trial points where the first constraint has not been
satisfied. The second line represents 11 points where the first
constraint has been satisfied but the second constraint has been
not. The last line shows 42 points where both constraints have
been satisfied and the objective function has been evaluated. The
total number of evaluations is equal to $21 + 11 \times 2 + 42
\times 3 = 169$. These evaluations have been executed during $21 +
11 + 42 = 74$ trials.

{}


\begin{thebibliography}{99}


\bibitem{Bertsekas (1996)}Bertsekas D.P. (1996), {\it Constrained Optimization and Lagrange Multiplier Methods}, Athena Scientific, Belmont, MA.

\bibitem{Famularo et al. (2002)}Famularo D., Sergeyev Ya.D., and Pugliese P. (2002), Test Problems  for Lipschitz Univariate Global Optimization with
Multiextremal Constraints. In: Dzemyda G., \v{S}altenis V., and
\v{Z}ilinskas A. (Eds.). {\it Stochastic and Global Optimization},
Kluwer Academic Publishers, Dordrecht, 93--110.

\bibitem{Horst and Pardalos (1995)}Horst R. and Pardalos P.M. (1995), {\it Handbook of Global Optimization}, Kluwer Academic Publishers, Dordrecht.

\bibitem{Nocedal and Wright (1999)}Nocedal J. and Wright S.J. (1999),  {\it Numerical Optimization} (Springer Series in Operations Research), Springer Verlag.

\bibitem{Pijavskii (1972)}Pijavskii S.A. (1972), An Algorithm for Finding the Absolute Extremum of a Function, {\it USSR Comput. Math. Math. Phys.}, {\bf 12} 57--67.

\bibitem{Pinter (1996)}Pint\'{e}r J.D. (1996), {\it Global Optimization in Action}, Kluwer Academic Publisher, Dordrecth.

\bibitem{Sergeyev (1995a)} Sergeyev Ya.D. (1995a) An information global optimization algorithm with local tuning, {\it SIAM J. Optim.} {\bf 5}, 858--870.

\bibitem{Sergeyev (1995b)} Sergeyev Ya.D. (1995b) A one-dimensional deterministic global minimization algorithm, {\it Comput. Math. Math. Phys.} {\bf 35}, 705--717.

\bibitem{Sergeyev (1998)}Sergeyev Ya.D. (1998), Global one-dimensional optimization using smooth auxiliary functions, {\it Math. Program.}, {\bf 81}, 127--146.

\bibitem{Sergeyev et al. (2001)}Sergeyev Ya.D., Famularo D., and
Pugliese P. (2001), Index Branch-and-Bound Algorithm for Lipschitz
Univariate Global Optimization with  Multiextremal Constraints,
{\it J. Global Optim.}, {\bf 21}, 317--341.

\bibitem{Sergeyev and Markin (1995)}Sergeyev Ya.D. and Markin D.L. (1995), An algorithm
for solving global optimization problems with nonlinear
constraints, {\it J. Global Optim.}, {\bf 7},  407--419.

\bibitem{Strongin (1978)}Strongin R.G. (1978), {\it Numerical Methods on Multiextremal Problems}, Nauka, Moscow. (In Russian).

\bibitem{Strongin (1984)}Strongin, R.G. (1984), Numerical methods for multiextremal nonlinear programming problems with nonconvex constraints. In: Demyanov V.F. and Pallaschke D. (Eds.). {\it Lecture Notes in Economics and Mathematical Systems} 255,  Proceedings 1984.  Springer-Verlag. IIASA,
Laxenburg/Austria, 278--282.

\bibitem{Strongin and Markin (1986)}Strongin R.G. and Markin D.L. (1986), Minimization of multiextremal functions with nonconvex constraints, {\it Cybernetics}, {\bf 22},  486--493.

\bibitem{Strongin and Sergeyev (2000)}Strongin R.G. and Sergeyev Ya.D. (2000), {\it Global Optimization with Non-Convex
Constraints: Sequential and Parallel Algorithms}, Kluwer Academic
Publishers, Dordrecht.

\end{thebibliography}
\end{document}